\documentclass[preprint,3p,sort]{elsarticle}

\journal{Journal of Computational Physics}

\usepackage{amsfonts}
\usepackage{amsthm}
\usepackage{graphicx}
\usepackage{epstopdf}
\usepackage{algorithmic}
\usepackage{amsmath}
\usepackage{amssymb}
\usepackage{hyperref}
\usepackage[T1]{fontenc}
\usepackage[dvipsnames]{xcolor}
\newcommand{\dbtilde}[1]{\widehat{#1}}
\newcommand{\BCMat}{\mathsf{R}}

\newtheorem{theorem}{Theorem}[section]
\newtheorem{lemma}[theorem]{Lemma}

\newtheorem{definition}[theorem]{Definition}
\theoremstyle{remark}
\newtheorem{remark}[theorem]{Remark}
\numberwithin{equation}{section}
\begin{document}

\begin{frontmatter}

%


\title{{A linear and nonlinear analysis of the shallow water equations and its impact on boundary conditions}}

\author[sweden,southafrica]{Jan Nordstr\"{o}m}
\author[sweden]{Andrew R.~Winters\corref{secondcorrespondingauthor}}
\cortext[secondcorrespondingauthor]{Corresponding author}
\ead{andrew.ross.winters@liu.se}
\address[sweden]{Department of Mathematics; Applied Mathematics, Link\"{o}ping University, SE-581 83 Link\"{o}ping, Sweden}
\address[southafrica]{Department of Mathematics and Applied Mathematics, University of Johannesburg, P.O. Box 524, Auckland Park 2006, Johannesburg, South Africa}

\begin{abstract}
We derive boundary conditions and estimates based on the energy and entropy analysis of systems of the nonlinear shallow water equations in two spatial dimensions. 
It is shown that the energy method provides more details, but is fully consistent with the entropy analysis. The details brought forward by the nonlinear energy analysis allow us to pinpoint where the difference between the linear and nonlinear analysis originate. We find that the result from the linear analysis does not necessarily hold in the nonlinear case. The nonlinear analysis leads in general to a different minimal number of boundary conditions compared with the linear analysis. In particular, and contrary to the linear case, the magnitude of the flow does not influence the number of required boundary conditions.
\end{abstract}



\end{frontmatter}

\section{Introduction}

Shallow water models provide a system of nonlinear hyperbolic equations that are useful to model fluid flows in lakes, rivers, oceans or near coastlines, e.g \cite{bonev2018discontinuous,ginting2019comparison,marras2018residual}. 
Specific information regarding the exact physical case is generally provided by boundary and initial conditions. 
Proper treatment of the boundary conditions is also an essential aspect of bounding the solution in the continuous analysis. This continuous analysis, in turn, is used as a guide to develop a stable numerical approximation for the shallow water model. For large scale problems, like those found in tsunami modelling {on a regional oceanic area}, the lateral boundaries of the domain are not physical boundaries. Therefore, artificial ``open''-type boundary conditions must be employed while maintaining proper solution estimates.

Due to the nonlinear nature of the shallow water equations (SWEs) the development of a bound of the solution becomes more subtle. In the literature one finds two general avenues in order to obtain estimates of solutions to nonlinear hyperbolic systems of conservation laws:
\begin{enumerate}
\item The energy method advocated by Kreiss \cite{kreiss1970,kreiss1989initial}, Gustafsson \cite{Gustafsson1978,gustafsson1995time} and Oliger and Sundstr\"{o}m \cite{oliger1978}.
\item The entropy stability theory used by Tadmor \cite{tadmor1984,Tadmor1987,Tadmor2003} and originating in works by Godunov \cite{godunov1961interesting}, Vol'pert \cite{volpert1967}, Kru\v{z}kov \cite{kruzkov1970}, Dafermos \cite{dafermos1973entropy}, Lax \cite{lax1973}, and Harten \cite{harten1983}.
\end{enumerate}

Traditionally, the energy method has been applied to linearized versions of systems of hyperbolic equations in order to develop boundary treatments that lead to energy estimates. In practice, these boundary conditions are needed to develop energy stable numerical approximations that weakly impose boundary information, e.g., through penalty terms \cite{carpenter1994time,nordstrom_roadmap,nordstrom2019,nordstrom2005} or numerical flux functions \cite{hindenlang2019,leveque1998,ESDGSEM2D_paper,xing2014}. The entropy method has often been applied to nonlinear hyperbolic systems on domains with periodic boundary conditions (or infinite domains) in order to obtain entropy conservation. This makes most investigations using entropy analysis similar to the classical von Neumann analysis in the sense that boundaries are often ignored \cite{Fjordholm2012_2,kuya2018kinetic,lefloch2000,merriam1989}. Entropy conservation is typically an intermediate step and dissipation is added to obtain entropy stability \cite{Tadmor2003}. The use of entropy stability theory to derive boundary conditions has been considered previously, see e.g.\cite{dalcin2019conservative,dubois1988,hindenlang2019,parsani2015entropy,svard2012,svard2014,svard2021entropy}, but mostly focused on solid walls. The main character of these two approaches are (roughly speaking): The energy analysis provides boundary conditions, but is difficult to apply in the nonlinear case. The entropy analysis is straightforward to apply for nonlinear equations, but its use for development of boundary conditions is limited.

In this work, we will derive stable nonlinear boundary treatments for the SWEs in two spatial dimensions. For this, we will examine the (physical) total energy of the system. This energy is an auxiliary conserved quantity, not explicitly built into the SWEs that can be used to create a nonlinear stability estimate \cite{gassner2015,tassi2007,ESDGSEM2D_paper} as well as build numerical approximations that respect the evolution of the total energy \cite{fjordholm2011,gassner2015,mishra2015}. An interesting quality of the SWEs and its nonlinear estimate is that the total energy also acts as a mathematical entropy function and fits into the entropy analysis framework \cite{Tadmor2003}.

Therefore, the total energy and the analysis of it for the SWEs act as a bridge between the classical energy method and the entropy analysis. So, a second contribution of this work is to examine the connection between the energy and entropy analysis, exemplifying their similarities and differences. We contrast the estimates from the nonlinear analysis to the ones from linear energy analysis and discuss their differences. In particular, we highlight that the linear analysis, in general, indicate a different number of boundary conditions than the nonlinear analysis. We also demonstrate that the nonlinear energy analysis provide additional details compared to the more easily applied entropy analysis, and makes it possible to pinpoint where the difference between the linear and nonlinear analysis reside.

We apply the energy method \cite{gustafsson1995time,kreiss1970,kreiss1989initial} and derive a bound for the total energy, which for the SWEs is a particular scaled version of the $L^2$ norm of the solution, based on a space and time-dependent weight matrix. It will be shown that the energy method is consistent with the mathematical entropy analysis, but also that it provides additional information, with respect to the boundary treatment. Investigations into energy stable boundary conditions for the linearized SWEs are many \cite{blayo2005,bristeau2001boundary,ghader2014,Mcdonald2002,oliger1978,shallowwaterbook} and leads to well-posed initial boundary value problems \cite{ghader2014,oliger1978,shallowwaterbook}. These linear boundary conditions are often applied in the nonlinear case \cite{ghader2014,verboom1984}. However, as we will show, the linear boundary treatments do not, in general, provide the correct minimal number of boundary conditions for the nonlinear shallow water problem.

The paper is organized as follows: Various forms of the SWEs are presented in Section~\ref{sec:sw}. {A brief review of energy estimates as well as linear or nonlinear boundedness and well-posedness is provided in Section~\ref{sec:semibound}.} An estimate of the total energy and entropy for the shallow water system is provided in Section~\ref{sec:energyStab}. In Section~\ref{sec:energyBCs}, we provide details, analysis, and discussion of the general open boundary conditions for the two-dimensional nonlinear SWEs in inflow and outflow regimes. In particular, we discuss similarities and differences between the results from the linear analysis and the nonlinear one. Concluding remarks are drawn in Section~\ref{sec:conclusion}.

\section{Shallow water equations}\label{sec:sw}

The SWEs over a flat bottom topography in conservative form \cite{whitham1974} are
\begin{equation}\label{eq:swCons}
\begin{aligned}
h_t + (hu)_x + (hv)_y &=0,\\[0.1cm]
(hu)_t + \left(hu^2 + \frac{g}{2}h^2\right)_{\! x} + (huv)_y - f h v &= 0,\\[0.1cm]
(hv)_t + (huv)_x + \left(hv^2 + \frac{g}{2}h^2\right)_{\! y} + f h u &= 0,
\end{aligned}
\end{equation}
which includes the continuity and momentum equations. Here $h(x,y,t)$ is the water height, $u(x,y,t)$ and $v(x,y,t)$ are the fluid velocities in the $x$- and $y$-directions, and $g$ is the gravitational constant. The system of equations \eqref{eq:swCons} are derived under the physical requirement that the water height is strictly positive, i.e., $h>0$ \cite{shallowwaterbook,whitham1974}. Additionally, we include the influence of Coriolis forces with the parameter $f$ which, for convenience, is assumed to be a constant. In practical applications $f$ is typically a function of latitude \cite{shallowwaterbook}, which would not affect the subsequent analysis in this work.

In order to apply the energy method it is convenient to work with the equivalent non-conservative form of the governing equations \cite{kreiss1977,nordstrom_roadmap,nordstrom2019,oliger1978} which is
\begin{equation}\label{eq:swNoncons}
    \begin{aligned}
    \phi_t + u\phi_x + v\phi_y + \phi (u_x + v_y) &=0,\\[0.1cm]
    u_t + u u_x + v u_y + \phi_x - f v &=0,\\[0.1cm]
    v_t + u v_x + v v_y + \phi_y + f u &=0.
    \end{aligned}
\end{equation}
In \eqref{eq:swNoncons}, we formulate the equations in terms of the geopotential $\phi = g h$ to simplify the analysis \cite{oliger1978}. Note, that $\phi>0$ according to the physical and mathematical requirements of the problem. Next, we write \eqref{eq:swNoncons} compactly in matrix-vector form by introducing the solution vector $\mathbf{q} = (\phi\,,\,u\,,\,v)^T$ as
\begin{equation*}
    \mathbf{q}_t + \mathcal{A}\mathbf{q}_x + \mathcal{B}\mathbf{q}_y+\mathcal{C}\mathbf{q} = \mathbf{0},
\end{equation*}
where
\begin{equation*}
    \mathcal{A} = \begin{bmatrix}
    u & \phi & 0\\[0.05cm]
    1 & u & 0\\[0.05cm]
    0 & 0 & u\\[0.05cm]
    \end{bmatrix},\qquad
    \mathcal{B} = \begin{bmatrix}
    v & 0 & \phi\\[0.05cm]
    0 & v & 0\\[0.05cm]
    1 & 0 & v\\[0.05cm]
    \end{bmatrix},\qquad
    \mathcal{C} = \begin{bmatrix}
    0 & 0 & 0\\[0.05cm]
    0 & 0 & f\\[0.05cm]
    0 & -f & 0\\[0.05cm]
    \end{bmatrix}.
\end{equation*}

The {specific} total energy is the sum of the kinetic and potential energy \cite{gassner2015}
\begin{equation}\label{eq:totEng}
\epsilon = \frac{\phi}{2 g}\left(u^2+v^2\right) + \frac{\phi^2}{2 g}.
\end{equation}
It is an auxiliary conserved quantity for smooth solutions of \eqref{eq:swCons} or \eqref{eq:swNoncons}. The {specific} total energy \eqref{eq:totEng} has associated flux functions \cite{fjordholm2011}
\begin{equation}\label{eq:entFluxes}
f^{\epsilon} = \frac{\phi u}{2g}\left(u^2+v^2\right) + \frac{\phi^2 u}{g},\qquad\textrm{and}\qquad g^{\epsilon} = \frac{\phi v}{2g}\left(u^2+v^2\right) + \frac{\phi^2 v}{g},
\end{equation}
that yield the  {specific} total energy conservation law
\begin{equation}\label{eq:energyConsLaw}
\epsilon_t + f^{\epsilon}_x + g^{\epsilon}_y = 0.
\end{equation}

In the entropy analysis, the  {specific} total energy in \eqref{eq:totEng} acts as a mathematical entropy function for the SWEs \cite{fjordholm2011}. As such, it is possible to define a new set of entropy variables $\mathbf{s}=(s_1,\,s_2,\,s_3)^T$ where
\begin{equation*}
s_1 = \frac{\partial \epsilon}{\partial h},\quad s_2 = \frac{\partial \epsilon}{\partial (hu)},\quad s_3 = \frac{\partial \epsilon}{\partial (hv)}.
\end{equation*}
Multiplying the conservative form of the SWEs \eqref{eq:swCons} from the left with $\mathbf{s}^T$ yields the auxiliary conservation law of the entropy function (in this case the {specific} total energy) \eqref{eq:energyConsLaw}. This contraction into entropy space involves the chain rule and certain compatibility conditions between the conservative fluxes from \eqref{eq:swCons} and the entropy fluxes from \eqref{eq:entFluxes} \cite{fjordholm2011,Tadmor2003}.

\begin{remark}
The Coriolis force is not present in the  {specific} total energy {conservation law} \eqref{eq:energyConsLaw}. This agrees with the underlying physics of the problem because the Coriolis terms do not perform work on the fluid \cite{shallowwaterbook}.
\end{remark}

{
\begin{remark}
The forthcoming nonlinear energy analysis considers the governing equations in non-conservative advective form \eqref{eq:swNoncons}. It is also possible to apply the nonlinear energy method to the conservative form \eqref{eq:swCons} although it more unwieldy in terms of algebraic manipulations required to obtain an energy statement at the boundary. Further, we will demonstrate in Section \ref{sec:energyStab} that the nonlinear energy analysis of the advective form is equivalent to the entropy analysis of the conservative equations as both arrive at the same final statement of the entropy flux in the normal direction at the boundary. However, we find that the nonlinear energy analysis reveals additional information to specifically see how each equation variable $(h, u, v)^T$ contributes to the boundary statement rather than the single expression that results form the entropy contraction. In essence, the entropy contraction ``hides'' an intermediate form of the boundary information that is revealed through the nonlinear energy analysis. This intermediate form is particularly useful to examine different open boundary configurations.
\end{remark}
}

{\section{Semi-boundedness and energy estimates for initial boundary value problems}\label{sec:semibound}}


{Before proceeding with the nonlinear energy analysis we review a few important concepts relevant for energy boundedness and well-posedness.}
{
Consider the generic initial boundary value problem (IBVP)
\begin{equation}
\begin{aligned}
	\mathbf{q}_t  + \mathsf{D}\mathbf{q} & = 0,
	& \quad & (x,y) \in \Omega, & \quad t & > 0,
	\\
	\mathsf{L} \mathbf{q} & =  0, & \quad & (x,y) \in \partial\Omega, & \quad t & > 0,
	\\
  	 \mathbf{q}  & =  f, & \quad & (x,y) \in \Omega, & \quad t & = 0
	\, .
\end{aligned}
\label{eq:semi_ibvp}
\end{equation}
In \eqref{eq:semi_ibvp}, $\mathbf{q} $ is the solution, $\mathsf{D}$ is a (possibly nonlinear) differential operator, $\mathsf{L}$ is the boundary operator, $\Omega$ is the domain, $\partial\Omega$ its boundary and $f$ is initial data.
}

{
For two vector functions $\mathbf{u}$ and $\mathbf{v}$ defined on $\Omega$, we introduce the inner product and corresponding norm
\begin{equation*}
	\langle \mathbf{u},\mathbf{v}\rangle = \int\limits_\Omega \mathbf{u}^T \mathcal{P} \mathbf{v} \,\mathrm{d}\Omega,
	\quad 
	\|\mathbf{u}\|^2 = \langle \mathbf{u},\mathbf{u}\rangle
	\, ,
\end{equation*}
where $ \mathcal{P}$ is a positive definite symmetric matrix. Further, we require the concept of semi-boundedness for linear problems \cite{kreiss1989initial,gustafsson2007,gustafsson1995time} and its extended definition to include nonlinear problems \cite{nordstrom2020spatial}. For linear problems we consider a space of functions and in the nonlinear case a set of functions.}  

{
\begin{definition}
\label{def:semi_bounded}
Let $\mathbb{V}$ be the space (or set) of all differentiable functions satisfying the boundary conditions $\mathsf{L}\mathbf{u} = 0$. The differential operator $\mathsf{D}$ in \eqref{eq:semi_ibvp} is {\normalfont semi-bounded} if for all $\mathbf{u}\in \mathbb{V}$ 
	\begin{equation}
		\langle \mathbf{u}, \mathsf{D}\mathbf{u} \rangle \ge \lambda \|\mathbf{u}\|^2
		\label{eq:semi_bounded}
	\end{equation}
holds, where $\lambda$ is a constant independent of $\mathbf{u}$.
\end{definition}
\begin{definition}
The differential operator $\mathsf{D}$ in \eqref{eq:semi_ibvp} is {\normalfont maximally semi-bounded} if it is semi-bounded in the function space (or function set) $\mathbb{V}$ but not in any space (set) with fewer boundary conditions.
\end{definition}
}

{
If $\mathsf{D}$ is maximally semi-bounded, we can estimate the solution to \eqref{eq:semi_ibvp} since 
\begin{equation}
\label{bounded_def}
	\| \mathbf{q} \|^2_t = 2\langle  \mathbf{q} , \mathbf{q} _t\rangle = -2 \langle  \mathbf{q} , \mathsf{D} \mathbf{q}  \rangle \le -2\lambda \| \mathbf{q} \|^2
	\Rightarrow
	\| \mathbf{q} (\cdot,T)\| \le e^{-\lambda T}\|f(\cdot)\|
	\, .
\end{equation}
\begin{remark}
The relation \eqref{eq:semi_bounded} for non-positive constants $\lambda$, yields a weak estimate. Homogeneous boundary conditions which lead to $\lambda > 0$ provide an estimate also for non-zero boundary data  \cite{nordstrom2019,nordstrom2020spatial}.
\end{remark}
}


{Next, we relate the concepts of semi-boundedness presented above to well-posedness which we roughly define as follows. The IBVP  \eqref{eq:semi_ibvp} is well posed if \textit{i)} the solution depends continuously on the given initial data (i.e. a bound like in (\ref{bounded_def}) exist), \textit{ii)} the solution is unique and \textit{iii)} the solution exists.}

{For linear IBVPs with smooth initial data and a maximally semi-bounded operator $\mathsf{D}$, well-posedness follows directly. Uniqueness follows by considering the difference between two solutions with the same data. The bound on the original problem bounds the difference and show that it is zero. Existence follows since no over-specification of boundary conditions is done (see \cite{nordstrom2020} for more details on how to determine the minimal number). Hence maximal semi-boundedness is a \textit{sufficient} condition for well-posedness of linear IBVPs. } 

{For nonlinear IBVPs we consider the analogous definition. We seek to determine the minimal number of boundary conditions such that the differential operator $\mathsf{D}$ is semi-bounded and lead to an energy estimate. As in the linear case, we then say that the operator is maximally semi-bounded. This is largely motivated by the linear theory, since the well-posedness theory for nonlinear hyperbolic systems is incomplete. Maximal semi-boundedness in the nonlinear case lead to an energy estimate, and hence it is a \textit{necessary} condition for well-posedness. However, since neither uniqueness nor existence follows, it is not a sufficient condition.
} 

\section{Energy and entropy analysis}\label{sec:energyStab}

Instead of using the energy (or entropy) conservation law \eqref{eq:energyConsLaw} directly, we will apply the classical energy method to the non-conservative system \eqref{eq:swNoncons}, which first require a suitable symmetrization matrix \cite{nordstrom_roadmap,oliger1978}. Guided by the linear analysis in \cite{oliger1978}, we select the symmetrizer
\begin{equation}\label{eq:symmMat}
    \mathcal{S} = \frac{1}{\sqrt{2g}}\begin{bmatrix}
    1 & 0 & 0\\[0.05cm]
    0 & \sqrt{\phi} & 0\\[0.05cm]
    0 & 0 & \sqrt{\phi}\\[0.05cm]
    \end{bmatrix}.
\end{equation}
The matrix $\mathcal{S}$ simultaneously symmetrizes the flux matrices $\mathcal{A}$ and $\mathcal{B}$ as
\begin{equation*}
    \mathcal{S}\mathcal{A}\mathcal{S}^{-1} = \begin{bmatrix}
    u & \sqrt{\phi} & 0\\[0.05cm]
    \sqrt{\phi} & u & 0\\[0.05cm]
    0 & 0 & u\\[0.05cm]
    \end{bmatrix},\qquad
    \mathcal{S}\mathcal{B}\mathcal{S}^{-1} = \begin{bmatrix}
    v & 0 & \sqrt{\phi}\\[0.05cm]
    0 & v & 0\\[0.05cm]
    \sqrt{\phi} & 0 & v\\[0.05cm]
    \end{bmatrix}.
\end{equation*}
Further, the symmetrization matrix directly relates the solution energy to the  {specific} total energy \eqref{eq:totEng} as
\begin{equation*}
    \left(\mathcal{S}\mathbf{q}\right)^T\mathcal{S}\mathbf{q} = \mathbf{q}^T\mathcal{P}\mathbf{q} = \frac{\phi u^2 + \phi v^2 + \phi^2}{2g},
\end{equation*}
where $\mathcal{P} = \mathcal{S}^2$.

Now we are equipped to apply the classical energy method \cite{gustafsson2007,gustafsson1995time,kreiss1970,nordstrom2019,nordstrom2005}. We pre-multiply \eqref{eq:swNoncons} by $\mathbf{q}^T\mathcal{P}$ to obtain
\begin{equation}\label{eq:contract1}
\mathbf{q}^T\mathcal{P}\mathbf{q}_t + \mathbf{q}^T\mathcal{P}\mathcal{A}\mathbf{q}_x + \mathbf{q}^T\mathcal{P}\mathcal{B}\mathbf{q}_y +\mathbf{q}^T\mathcal{P}\mathcal{C}\mathbf{q} = 0.
\end{equation}
From the skew-symmetry of the Coriolis matrix we immediately see that $\mathbf{q}^T\mathcal{P}\mathcal{C}\mathbf{q} = 0$.
The flux matrices are now symmetrized and take the form
\begin{equation}\label{eq:symmetrizeSWE}
    \mathcal{P}\mathcal{A} = \frac{1}{2g}\begin{bmatrix}
    u & \phi & 0\\[0.05cm]
    \phi & \phi u & 0\\[0.05cm]
    0 & 0 & \phi u\\[0.05cm]
    \end{bmatrix},\qquad
    \mathcal{P}\mathcal{B} = \frac{1}{2g}\begin{bmatrix}
    v & 0 & \phi\\[0.05cm]
    0 & \phi v & 0\\[0.05cm]
    \phi & 0 & \phi v\\[0.05cm]
    \end{bmatrix}.
\end{equation}
We seek to rewrite \eqref{eq:contract1} with complete derivatives, and use the relations
\begin{equation*}
\begin{aligned}
    \mathbf{q}^T\mathcal{P}\mathbf{q}_t &= \frac{1}{2}\left(\mathbf{q}^T\mathcal{P}\mathbf{q}\right)_{\! t} -\frac{1}{2}\mathbf{q}^T\mathcal{P}_t\mathbf{q},\\[0.1cm]
    \mathbf{q}^T\mathcal{P}{\mathcal{A}}\mathbf{q}_x &= \frac{1}{2}\left(\mathbf{q}^T\mathcal{P}{\mathcal{A}}\mathbf{q}\right)_{\! x} -\frac{1}{2}\mathbf{q}^T\left(\mathcal{P}{\mathcal{A}}\right)_x\mathbf{q},\\[0.1cm]
    \mathbf{q}^T\mathcal{P}{\mathcal{B}}\mathbf{q}_y &= \frac{1}{2}\left(\mathbf{q}^T\mathcal{P}{\mathcal{B}}\mathbf{q}\right)_{\! y} -\frac{1}{2}\mathbf{q}^T\left(\mathcal{P}{\mathcal{B}}\right)_y\mathbf{q}.
    \end{aligned}
\end{equation*}
The expression \eqref{eq:contract1} becomes
\begin{equation}\label{eq:contract2}
    \left(\mathbf{q}^T\mathcal{P}\mathbf{q}\right)_{\! t} + \left(\mathbf{q}^T\mathcal{P}{\mathcal{A}}\mathbf{q}\right)_{\! x} + \left(\mathbf{q}^T\mathcal{P}{\mathcal{B}}\mathbf{q}\right)_{\! y} - \mathbf{q}^T\left(\mathcal{P}_t + \left(\mathcal{P}{\mathcal{A}}\right)_x + \left(\mathcal{P}{\mathcal{B}}\right)_y\right)\mathbf{q} = 0.
\end{equation}
%

\subsection{Nonlinear energy analysis}
%

To begin, we compute the derivatives of the matrices in \eqref{eq:contract2}  to be
\begin{equation*}
\begin{aligned}
    \mathcal{P}_t = \frac{1}{2g}\begin{bmatrix}
    0 & 0 & 0\\[0.05cm]
    0 & \phi_t & 0\\[0.05cm]
    0 & 0 & \phi_t\\[0.05cm]
    \end{bmatrix},\quad
    \left(\mathcal{P}{\mathcal{A}}\right)_x &= \frac{1}{2g}\begin{bmatrix}
    u_x & \phi_x & 0\\[0.05cm]
    \phi_x & \phi_x u + \phi u_x & 0\\[0.05cm]
    0 & 0 & \phi_x u + \phi u_x\\[0.05cm]
    \end{bmatrix},\\[0.1cm]
    \left(\mathcal{P}{\mathcal{B}}\right)_y &= \frac{1}{2g}\begin{bmatrix}
    v_y & 0 & \phi_y \\[0.05cm]
    0 & \phi_y v + \phi v_y & 0\\[0.05cm]
    \phi_y & 0 & \phi_y v + \phi v_y\\[0.05cm]
    \end{bmatrix}
\end{aligned}
\end{equation*}
which gives, noting the continuity equation (the first equation) from \eqref{eq:swNoncons},
\begin{equation*}
    \mathcal{P}_t+\left(\mathcal{P}{\mathcal{A}}\right)_x+\left(\mathcal{P}{\mathcal{B}}\right)_y = \frac{1}{2g}\begin{bmatrix}
    u_x+v_y & \phi_x & \phi_y\\[0.05cm]
    \phi_x & 0 & 0\\[0.05cm]
    \phi_y & 0 & 0\\[0.05cm]
    \end{bmatrix}.
\end{equation*}

It is now straightforward to compute
\begin{equation*}
\begin{aligned}
    \mathbf{q}^T\left(\mathcal{P}_t+\left(\mathcal{P}{\mathcal{A}}\right)_x+\left(\mathcal{P}{\mathcal{B}}\right)_y\right)\mathbf{q} &= \frac{1}{2g}\left[\phi^2\left(u_x + v_y\right) + 2\phi\phi_x u+ 2\phi\phi_y v\right]\\[0.05cm]
    &=\frac{1}{2g}\left[\left(\phi^2 u\right)_{\! x} +\left(\phi^2 v\right)_{\! y}\right].
    \end{aligned}
\end{equation*}
Therefore, \eqref{eq:contract2} becomes
\begin{equation}\label{eq:energyEqn1}
    \left(\mathbf{q}^T\mathcal{P}\mathbf{q}\right)_{\! t} + \left(\mathbf{q}^T\mathcal{P}{\mathcal{A}}\mathbf{q}\right)_{\! x} + \left(\mathbf{q}^T\mathcal{P}{\mathcal{B}}\mathbf{q}\right)_{\! y} -\frac{1}{2g}\left[\left(\phi^2 u\right)_{\! x} +\left(\phi^2 v\right)_{\! y}\right] = 0.
\end{equation}

We continue by introducing the norm $\|\mathbf{q}\|^2_{\mathcal{P}} = \int_{\Omega}\mathbf{q}^T\mathcal{P}\mathbf{q}\,\mathrm{d}x\mathrm{d}y$, integrating \eqref{eq:energyEqn1} over $\Omega\subset\mathbb{R}^2$, and apply Gauss' theorem to find
\begin{equation}\label{eq:boundaryPart1}
\frac{d}{dt}\|\mathbf{q}\|^2_{\mathcal{P}} + \frac{1}{2g}\oint\limits_{\partial\Omega}\left(\mathbf{q}^T\widetilde{\mathcal{A}}\mathbf{q} -  \frac{\phi^2}{2g}u_n\right)\,\mathrm{d}S=0,
\end{equation}
where $u_n = u n_x + v n_y$ is the normal velocity to the boundary $\partial\Omega$ and
\begin{equation*}
    \frac{1}{2g}\widetilde{\mathcal{A}} = \mathcal{P}\left(\mathcal{A}n_x + \mathcal{B}n_y\right) = \frac{1}{2g}\begin{bmatrix}
    u_n & n_x \phi & n_y \phi\\[0.05cm]
    n_x \phi & \phi u_n & 0\\[0.05cm]
    n_y \phi & 0 & \phi u_n\\[0.05cm]
    \end{bmatrix}.
\end{equation*}

{To simplify the upcoming analysis, we rotate the solution vector $\mathbf{q}$ as well as the matrix $\widetilde{\mathcal{A}}$ into the normal and tangential velocity directions, using
\begin{equation*}
\mathcal{T} = \begin{bmatrix}
1 & 0 & 0\\[0.05cm]
0 & n_x & n_y\\[0.05cm]
0 & -n_y & n_x\\[0.05cm]
\end{bmatrix},\quad\mathrm{with}\quad \mathcal{T}^{-1} = \mathcal{T}^T,
\end{equation*}
such that
\begin{equation*}
\mathcal{T}\mathbf{q} = \begin{bmatrix}
\phi\\[0.05cm]
n_x u + n_y v\\[0.05cm]
-n_y u + n_x v\\[0.05cm]
\end{bmatrix}
= \begin{bmatrix}
\phi\\[0.05cm]
u_n\\[0.05cm]
u_s\\[0.05cm]
\end{bmatrix}
= \widetilde{\mathbf{q}},
\end{equation*}
and
\begin{equation*}
\mathbf{q}^T\widetilde{\mathcal{A}}\mathbf{q} = \mathbf{q}^T\mathcal{T}^T\left(\mathcal{T}\widetilde{\mathcal{A}}\mathcal{T}^T\right)\mathcal{T}\mathbf{q} =  \widetilde{\mathbf{q}}^T\dbtilde{\mathcal{A}}\widetilde{\mathbf{q}}\,,
\qquad
\dbtilde{\mathcal{A}}
=
\begin{bmatrix}
u_n & \phi & 0\\[0.05cm]
\phi & \phi u_n & 0\\[0.05cm]
0 & 0 &\phi u_n \\[0.05cm]
\end{bmatrix}.
\end{equation*}
The energy evolution equation \eqref{eq:boundaryPart1} can now be rewritten as 
\begin{equation}\label{eq:boundaryPart2}
\frac{d}{dt}\|\mathbf{q}\|^2_{\mathcal{P}} + \frac{1}{2g}\oint\limits_{\partial\Omega}\left(\widetilde{\mathbf{q}}^T\dbtilde{\mathcal{A}}\widetilde{\mathbf{q}} -  \frac{\phi^2}{2g}u_n\right)\,\mathrm{d}S=0.
\end{equation}
}

{Next, we examine the effect of the additional scalar term in \eqref{eq:boundaryPart2} introduced by the nonlinearity of the problem. 
\begin{lemma}\label{lemma:Nmatrix}
There are three unique ways to rewrite the scalar term in \eqref{eq:boundaryPart2} as a quadratic form
\begin{equation}\label{eq:nonlinearTerminLemma}
 \frac{\phi^2}{2g}u_n = \frac{1}{2g}\widetilde{\mathbf{q}}^T\mathcal{N}_i\widetilde{\mathbf{q}},\qquad i = 1,2,\,\mathrm{or}\;3
\end{equation}
where
\begin{equation}\label{eq:threeMatrices}
\mathcal{N}_1
=
\begin{bmatrix}
u_n & 0 & 0\\[0.05cm]
0 & 0 & 0\\[0.05cm]
0 & 0 & 0\\[0.05cm]
\end{bmatrix},
\;\;
\mathcal{N}_2
=
\begin{bmatrix}
0 & \frac{\phi}{2} & 0\\[0.05cm]
\frac{\phi}{2} & 0 & 0\\[0.05cm]
0 & 0 & 0\\[0.05cm]
\end{bmatrix}
\;\;\mathrm{ and }\;\;
\mathcal{N}_3
=
\begin{bmatrix}
0 & 0 & 0\\[0.05cm]
0 & \frac{\phi^2}{u_n} & 0\\[0.05cm]
0 & 0 & 0\\[0.05cm]
\end{bmatrix}.
\end{equation}
\end{lemma}
\begin{proof}
See \ref{sec:newMatrix}.
\end{proof}
}

{The results of Lemma \ref{lemma:Nmatrix} allows the reformulation of \eqref{eq:boundaryPart2} to become
\begin{equation}\label{eq:boundaryPart3}
\frac{d}{dt}\|\mathbf{q}\|^2_{\mathcal{P}} + \frac{1}{2g}\oint\limits_{\partial\Omega}\widetilde{\mathbf{q}}^T(\dbtilde{\mathcal{A}}- \mathcal{N}_i)\widetilde{\mathbf{q}}\,\mathrm{d}S=0,\quad i = 1,2,3.
\end{equation}
}

\begin{remark}
If we use the linear analysis and neglect the matrices $\mathcal{N}_i$, the contribution to the energy will be different, which leads to other boundary procedures \cite{oliger1978,shallowwaterbook} (that will be discussed in detail in Section~\ref{sec:linearenergyStab}).
\end{remark}

To develop an energy estimate we must bound the surface integral term in \eqref{eq:boundaryPart3}. To find such a bound we rotate the boundary matrix $\dbtilde{\mathcal{A}}- \mathcal{N}_i$ into a diagonal form \cite{nordstrom_roadmap,nordstrom2019,nordstrom2020,nordstrom2005}. 
For each available boundary matrix, we introduce a solution dependent rotation matrix and corresponding vector of linearly independent rotated variables, i.e.,
\begin{equation}\label{eq:generalRotation}
\mathcal{R}_i = \begin{bmatrix}
\alpha_i & 0 & 0\\[0.05cm]
\beta_i & \gamma_i & 0\\[0.05cm]
0 & 0 & \delta_i\\[0.05cm]
\end{bmatrix}
\quad
\text{ and }
\quad
\mathbf{w}_i = 
\mathcal{R}_i^{-1}\widetilde{\mathbf{q}}
=
\begin{bmatrix}
\frac{\phi}{\alpha_i}\\[0.05cm]
\frac{1}{\gamma_i}\left(u_n - \frac{\beta_i}{\alpha_i}\phi\right)\\[0.05cm]
\frac{u_s}{\delta_i}\\[0.05cm]
\end{bmatrix}.
\end{equation}
Each rotation matrix $\mathcal{R}_i$, $i = 1,2,3$ transform the boundary matrix into a diagonal form
\begin{equation*}
\mathcal{R}_i^T(\dbtilde{\mathcal{A}}- \mathcal{N}_i)\mathcal{R}_i = \mathcal{D}_i,
\end{equation*}
such that the boundary integral term in \eqref{eq:boundaryPart3} becomes
\begin{equation*}
\frac{1}{2g}
\oint\limits_{\partial\Omega}\widetilde{\mathbf{q}}^T(\dbtilde{\mathcal{A}}- \mathcal{N}_i)\widetilde{\mathbf{q}}\,\mathrm{d}S 
=
\frac{1}{2g}
\oint\limits_{\partial\Omega}\mathbf{w}_i^T\mathcal{D}_i\mathbf{w}_i\,\mathrm{d}S,\quad i = 1,2,3.
\end{equation*}
We collect the rotations of the form \eqref{eq:generalRotation} for the three boundary matrices:
\begin{description}
\item[{\fboxsep=1pt\fbox{$\dbtilde{\mathcal{A}} - \mathcal{N}_1$}}] The terms in $\mathcal{R}_1$ are $\alpha_1 = 1$, $\beta_1 = -1/u_n$ and $\gamma_1 = \delta_1 = 1/u_n$. The rotated variables are $\mathbf{w}_1 = \left[\phi,u_n^2 + \phi,u_n u_s\right]^T$ with the diagonal matrix
\begin{equation}\label{eq:DN1}
\mathcal{D}_1 = \frac{\phi}{u_n}\text{diag}\left(-1,1,1\right).
\end{equation}
\item[{\fboxsep=1pt\fbox{$\dbtilde{\mathcal{A}} - \mathcal{N}_2$}}] For $\mathcal{R}_2$, the values $\alpha_2 = -2$ and $\beta_2 = \gamma_2 = \delta_2 = 1/u_n$ in \eqref{eq:generalRotation} give the rotated variables $\mathbf{w}_2 = \left[-\phi/2,u_n^2 + \phi/2,u_n u_s\right]^T$ and diagonal matrix
\begin{equation}\label{eq:DN2}
\mathcal{D}_2 = \frac{\phi}{u_n}\text{diag}\left(\frac{4u_n^2 - \phi}{\phi},1,1\right).
\end{equation}
\item[{\fboxsep=1pt\fbox{$\dbtilde{\mathcal{A}} - \mathcal{N}_3$}}] The third matrix $\mathcal{R}_3$ takes the values
\begin{equation*}
\alpha_3 = 1,\quad \beta_3 = -\frac{u_n}{u_n^2-\phi},\quad \gamma_3 = \frac{u_n}{u_n^2-\phi} 
,\quad \delta_3 = \frac{1}{u_n},
\end{equation*}
which yields the rotated variables $\mathbf{w}_3 = \left[\phi,u_n^2,u_n u_s\right]^T$ and diagonal matrix
\begin{equation}\label{eq:DN3}
\mathcal{D}_3 = \frac{\phi}{u_n}\text{diag}\left(\frac{u_n^2(u_n^2-2\phi)}{\phi(u_n^2-\phi)},\frac{u_n^2}{u_n^2-\phi},1\right).
\end{equation}
\end{description}

\begin{remark}
The contracted, rotated forms $\mathbf{w}_i^T\mathcal{D}_i\mathbf{w}_i$, $i = 1,2,3$ are equivalent. However, only the boundary matrix associated with $\mathcal{N}_1$ can be rotated such that one can, by inspection, determine the number of positive and negative diagonal values. The reason is that the funtional dependence of the diagonal values in \eqref{eq:DN1} can be factored out, clearly displaying the signs of the diagonal values. In contrast, the entries in \eqref{eq:DN2} and \eqref{eq:DN3} are complicated functions of the solution variables $\phi$ and $u_n$, that cannot be factored out. Therefore, determining the number of positive and negative entries in either $\mathcal{D}_2$ or $\mathcal{D}_3$ requires careful examination of the size and sign of the normal velocity in comparison to the geopotential. In \eqref{eq:DN2} and \eqref{eq:DN3} one must consider the diagonal values and the rotated variables at \textit{the same time}. One cannot, as in \eqref{eq:DN1}, \textit{first} determine the sign of the diagonal values and \textit{next} combine the variables to create boundary conditions.
\end{remark}

For the upcoming analysis of stable boundary conditions in Section~\ref{sec:energyBCs} we select the first rotation to have \eqref{eq:DN1}, and \eqref{eq:boundaryPart3} becomes
\begin{equation}\label{eq:boundaryPartFinal2}
\frac{d}{dt}\|\mathbf{q}\|^2_{\mathcal{P}} + \frac{1}{2 g}\oint\limits_{\partial\Omega}\mathbf{w}_1^T\mathcal{D}_1\mathbf{w}_1\,\mathrm{d}S=0.
\end{equation}

\subsection{Nonlinear entropy analysis}

It is interesting to examine how the statement \eqref{eq:energyEqn1} compares to the energy and entropy conservation law \eqref{eq:energyConsLaw}. Due to the construction of the symmetrization matrix \eqref{eq:symmMat}, we have the time evolution of the  {specific} total energy
\begin{equation*}
    \left(\mathbf{q}^T\mathcal{P}\mathbf{q}\right)_{\! t} = \left(\frac{\phi u^2 + \phi v^2 + \phi^2}{2g}\right)_{\! t} = \epsilon_t.
\end{equation*}
We look at the energy flux contribution in the $x$-direction to find
\begin{equation*}
\begin{aligned}
    \left(\mathbf{q}^T\mathcal{P}{\mathcal{A}}\mathbf{q}\right)_{\! x} - \frac{1}{2g}\left(\phi^2 u\right)_{\! x} &= \left(\frac{\phi u^3}{2g} + \frac{\phi u v^2}{2g} + \frac{\phi^2 u}{g}+ \frac{\phi^2 u}{2g}\right)_{\! x} - \frac{1}{2g}\left(\phi^2 u\right)_{\! x}\\[0.05cm]
    &=\left(\frac{\phi u}{2g}\left(u^2+v^2\right) + \frac{\phi^2 u}{g}\right)_{\! x}=f^{\epsilon}_x.
    \end{aligned}
\end{equation*}
Similarly, in the $y$-direction we have
\begin{equation*}
    \left(\mathbf{q}^T\mathcal{P}{\mathcal{B}}\mathbf{q}\right)_{\! y} - \frac{1}{2g}\left(\phi^2 v\right)_{\! y} =\left(\frac{\phi v}{2g}\left(u^2+v^2\right) + \frac{\phi^2 v}{g}\right)_{\! y} = g^{\epsilon}_y.
\end{equation*}
So, we find that \eqref{eq:energyEqn1} is equivalent to the entropy conservation law \eqref{eq:energyConsLaw}. Integrating \eqref{eq:energyConsLaw} over a general domain and applying Gauss' theorem \cite{evans2010} yield
\begin{equation}\label{eq:energyEq}
\int\limits_{\Omega}\left(\epsilon_t + f^{\epsilon}_x + g^{\epsilon}_y\right)\,\mathrm{d}\mathbf{x} = \int\limits_{\Omega} \epsilon_t \,\mathrm{d}\mathbf{x} + \oint\limits_{\partial\Omega} f^{\epsilon} n_x + g^{\epsilon} n_y\,\mathrm{d}S = 0,
\end{equation}
with the outward pointing normal vector $\vec{n}=(n_x,n_y)^T$. The energy conservation law \eqref{eq:energyEq} is the same as \eqref{eq:boundaryPartFinal2} albeit derived from a different perspective. We expand the boundary term to highlight the similarity to \eqref{eq:boundaryPartFinal2}
\begin{equation*}
\int\limits_{\Omega} \epsilon_t \,\mathrm{d}\mathbf{x} + \frac{1}{2g}\oint\limits_{\partial\Omega}\phi u_n(2\phi+u_n^2+u_s^2)\,\mathrm{d}S = 0.
\end{equation*}

\begin{remark}
For non-smooth solutions, the  {specific} total energy conservation law \eqref{eq:energyEq} becomes the energy (or entropy) inequality
\begin{equation*}
\int\limits_{\Omega} \epsilon_t + \frac{1}{2g}\oint\limits_{\partial\Omega}\phi u_n(2\phi+u_n^2+u_s^2)\,\mathrm{d}S \leq 0,
\end{equation*}
which states that the  {specific} total energy (or entropy) dissipates in the presence of discontinuities \cite{fjordholm2011,gassner2015}. It is the typical entropy condition for conservation laws \cite{dafermos2005hyperbolic,toro2013riemann}.
\end{remark}

\begin{remark}
The analysis above shows that the nonlinear energy and entropy analysis are similar, but that the energy analysis provides more details. In particular, it provides the matrix $\mathcal{P}$ forming the norm, the symmetrized matrices $\mathcal{P}\mathcal{A}$ and $\mathcal{P}\mathcal{B}$ in \eqref{eq:symmetrizeSWE} (that depend on space and time but have the same form as in the linear analysis \cite{oliger1978}) and an explicit matrix-vector formulation of the boundary term. The energy analysis also shows that the difference between the linear and nonlinear analysis depends only on $\mathcal{N}_i$ in \eqref{eq:nonlinearTerminLemma}. The entropy analysis on the other hand is more direct and easily applied. Having established the similarities and differences between the energy and entropy analysis, we in the following refer to both when we state energy (or entropy).
\end{remark}

\subsection{Linear energy analysis}\label{sec:linearenergyStab}

The standard linear energy analysis of the SWE is \textit{contained} within equation \eqref{eq:contract2}. This is because the nonlinearity of the problem produces the last term in \eqref{eq:contract2} which originates from the space- and time-dependent symmetrizer $\mathcal{S}$ in \eqref{eq:symmMat}. If one linearizes the equations \eqref{eq:swNoncons} around a constant background state, the coefficient matrices $\mathcal{A}$ and $\mathcal{B}$ as well as the symmetrizer $\mathcal{S}$ become constants with the same functional form as in the nonlinear case and the final term in \eqref{eq:contract2} vanishes. 

Ignoring the nonlinear term $\mathcal{N}_i$ and performing such a linearization changes the form of the boundary contribution from \eqref{eq:boundaryPart3} to
\begin{equation}\label{eq:boundaryPartLinear}
\frac{d}{dt}\|\mathbf{q}\|^2_{\mathcal{P}} + \frac{1}{2g}\oint\limits_{\partial\Omega}\left(\bar{u}_n\phi^2 +  \overline{\phi}\,\widetilde{\mathbf{q}}^T
\begin{bmatrix}
0 & 1 & 0\\[0.05cm]
1 &\bar{u}_n & 0\\[0.05cm]
0 & 0 & \bar{u}_n\\[0.05cm]
\end{bmatrix}
\widetilde{\mathbf{q}}\right)\mathrm{d}S = 0,
\end{equation}
where $\overline{\phi}$ and $\bar{u}_n$ are constant mean states for the geopotential and normal velocity. We use a rotation like \eqref{eq:generalRotation} to bring the linearized matrix term in \eqref{eq:boundaryPartLinear} into a diagonal form
\begin{equation*}
\overline{\mathcal{D}}
= 
\overline{\mathcal{R}}^{T}
\begin{bmatrix}
0 & 1 & 0\\[0.05cm]
1 & \bar{u}_n & 0\\[0.05cm]
0 & 0 & \bar{u}_n\\[0.05cm]
\end{bmatrix}
\overline{\mathcal{R}}
=
\text{diag}\left(-\frac{1}{\bar{u}_n},\frac{1}{\bar{u}_n},\frac{1}{\bar{u}_n}\right)
\;\;\text{ with }\;\;
\overline{\mathcal{R}}
= 
\begin{bmatrix}
1 & 0 & 0\\[0.05cm]
-\frac{1}{\bar{u}_n} & \frac{1}{\bar{u}_n} & 0\\[0.05cm]
0 & 0 & \frac{1}{\bar{u}_n}\\[0.05cm]
\end{bmatrix}.
\end{equation*}
This also gives the linearized rotated variables $\overline{\mathbf{w}} = \overline{\mathcal{R}}^{-1}\widetilde{\mathbf{q}} = [\phi\,,\,\bar{u}_nu_n + \phi\,,\,\bar{u}_n u_s]^T$. Now the boundary term can be compactly written as
\begin{equation}\label{eq:boundaryPartLinear2}
\frac{d}{dt}\|\mathbf{q}\|^2_{\mathcal{P}} + \frac{1}{2g}\overline{\phi} \oint\limits_{\partial\Omega}\left(\bar{u}_n\phi^2 + \overline{\phi}\,\overline{\mathbf{w}}^T\overline{\mathcal{D}}\,\overline{\mathbf{w}}\right)\,\mathrm{d}S=0.
\end{equation}
Next, we manipulate the first term of the boundary integral in \eqref{eq:boundaryPartLinear2} to get
\begin{equation}\label{eq:diagonalLinearMat}
\bar{u}_n\phi^2 + \overline{\phi}\,\overline{\mathbf{w}}^T\overline{\mathcal{D}}\,\overline{\mathbf{w}} 
= \frac{\overline{\phi}}{\bar{u}_n}\frac{\bar{u}^2_n}{\overline{\phi}}\phi^2 + \overline{\phi}\,\overline{\mathbf{w}}^T\overline{\mathcal{D}}\,\overline{\mathbf{w}}
=
\frac{\overline{\phi}}{\bar{u}_n}
\overline{\mathbf{w}}^T
\begin{bmatrix}
\frac{\bar{u}^2_n - \overline{\phi}}{\overline{\phi}} & 0 & 0\\[0.05cm]
0 & 1 & 0\\[0.05cm]
0 & 0 & 1\\[0.05cm]
\end{bmatrix}
\overline{\mathbf{w}}.
\end{equation}

Now the sign of the diagonal matrix 
shifts with the sign of $\bar{u}_n^2 - \overline{\phi}$. When $\bar{u}_n^2 < \overline{\phi}$ the flow is subcritical whereas when $\bar{u}_n^2 > \overline{\phi}$, we have supercritical flow \cite{shallowwaterbook}. {The matrix in \eqref{eq:diagonalLinearMat} is scaled by the value $\overline{\phi}/\bar{u}_n$ which further influences the overall sign of each term.} So we see that the signs of the matrix in \eqref{eq:diagonalLinearMat} progresses through four states depending on the direction and magnitude of the normal velocity \cite{huang2014linearized}:
\begin{itemize}
\renewcommand\labelitemi{\textbullet}
\item Supercritical inflow where $\bar{u}_n < 0$ yields three negative values.
\item Subcritical inflow where $\bar{u}_n < 0$ yields two negative values.
\item Supercritical outflow where $\bar{u}_n > 0$ yields one negative value.
\item Subcritical outflow where $\bar{u}_n > 0$ yields zero negative values.
\end{itemize}
The impact such sign changes have on the number of boundary conditions is discussed in the next Section. Linear boundary treatments for the SWE were originally developed by Oliger and Sundstr\"{o}m \cite{oliger1978} and expanded upon elsewhere, e.g., by Ghader and Nordstr\"{o}m \cite{ghader2014} and Huang and Temam \cite{huang2014linearized}. 

\begin{remark}
Using lifting approaches, e.g. \cite{nordstrom_roadmap}, it is straightforward to apply the results from the forthcoming analysis and implement weak boundary conditions for numerical calculations \cite{carpenter1994time,nordstrom_roadmap,nordstrom2019,nordstrom2005,hindenlang2019,leveque1998,ESDGSEM2D_paper,xing2014}.
\end{remark}

\section{Energy and entropy stable boundary conditions}\label{sec:energyBCs}

The sign of the normal velocity determines whether there are inflow or outflow conditions at the domain boundary. That is, $u_n<0$ corresponds to inflow conditions, $u_n > 0$ to outflow conditions and $u_n=0$ to {glancing conditions (whose prototypical physical example is a rigid boundary)}.
Energy stable open boundary conditions of Dirichlet type for a general domain must be chosen such that the reformulated boundary term in \eqref{eq:DN1} and \eqref{eq:boundaryPartFinal2} 
\begin{equation}\label{eq:rewriteBCs}
\mathbf{w}^T\mathcal{D}\,\mathbf{w}
=
\begin{bmatrix}
\mathbf{w}^{+}\\
\mathbf{w}^{-}\\
\end{bmatrix}^{\!T}
\!
\left[\begin{array}{@{}c|c@{}}
\mathcal{D}^{+} & 0\\\hline
0 & \mathcal{D}^{-}\\
\end{array}\right]
\!
\begin{bmatrix}
\mathbf{w}^{+}\\
\mathbf{w}^{-}\\
\end{bmatrix}
=
\frac{\phi}{u_n}
\begin{bmatrix}
\phi \\
u_n^2 + \phi\\
u_n u_s\\
\end{bmatrix}^{\!T}
\!
\left[\begin{array}{@{}ccc@{}}
-1& 0 & 0\\[0.1cm]
0 & 1 & 0\\[0.1cm]
0 & 0 & 1 \\[0.1cm]
\end{array}\right]
\!
\begin{bmatrix}
\phi \\
u_n^2 + \phi\\
u_n u_s\\
\end{bmatrix}
\end{equation}
is made positive semi-definite. 
{That is, a minimum number of boundary conditions are determined so that 
\begin{equation} 
\mathbf{w}^T\mathcal{D}\,\mathbf{w} \geq 0
\label{eq:rewriteBCs_inequality}
\end{equation}
holds, which means that the operator $\mathcal{D}$ is maximally semi-bounded.
}

Above we know that the geopotential $\phi$ is strictly positive by the physical requirements, and we have supressed the subscript $1$. 
In \eqref{eq:rewriteBCs}, $\mathbf{w}^{+}$ contains the outgoing boundary information and $\mathbf{w}^{-}$ the incoming information. The diagonal blocks {$\mathcal{D}^{+}$ and $\mathcal{D}^{-}$ contain the positive and negative diagonal entries}, respectively. Note that the signs of the entries in the matrix $\mathcal{D}$ \eqref{eq:DN1} depend entirely on the sign of the normal velocity. This clear separation of the positive and negative contributions {to the energy at the boundary} provides a suitable setting (not present in \eqref{eq:DN2} or \eqref{eq:DN3}) for discussing boundary conditions for the nonlinear problem. It will also serve to illustrate the difference between a linear and nonlinear analysis.

Our goal is to guarantee energy stability at open boundaries with a minimal number of boundary conditions. As a starting point, we follow the linear analysis path and assume that the number of boundary conditions required is equal to the size of $\mathcal{D}^{-}$ \cite{nordstrom_roadmap,nordstrom2020}. {However}, in the course of the forthcoming analysis we will demonstrate that this assumption can be misleading. That is, it may be the case that more (or fewer) boundary conditions are required to guarantee that the inequality {\eqref{eq:rewriteBCs_inequality}} holds. The subtlety of the nonlinear analysis is that the entries of the matrix $\mathcal{D}$ in general \textit{depend} upon the solution quantities (for instance the normal velocity)
, whereas such a dependency does not exist in the linear case \cite{nordstrom_roadmap}.
\begin{remark}
As discussed in Section \ref{sec:semibound}, an energy estimate can always be obtained by specifying too many boundary conditions \cite{nordstrom_roadmap,nordstrom2020}. However, in general that means that existence cannot be obtained for linear problems. The precise existence requirement is not known for nonlinear problems, but it is reasonable to assume (unless proven otherwise) that similar requirements as for linear problems hold.
\end{remark}

A general form of the boundary conditions is \cite{nordstrom2019}
\begin{equation}\label{eq:generalFormBCs}
\mathbf{w}^{-} = \BCMat\mathbf{w}^{+} + \mathbf{g},
\end{equation}
where $\BCMat$ is a coefficient matrix with the number of rows equal to the minimal number of required boundary conditions. The vector $\mathbf{g}$ contains known external data at the boundary. Essentially, the boundary condition \eqref{eq:generalFormBCs} represents incoming boundary information as a linear combination of the outgoing information and external data. With the rewritten stability condition \eqref{eq:rewriteBCs}, {its corresponding inequality \eqref{eq:rewriteBCs_inequality}} and the general boundary condition \eqref{eq:generalFormBCs} we have
\begin{theorem}\label{thm:nonLinStabBCs}
The inhomogeneous, general boundary conditions \eqref{eq:generalFormBCs} for the nonlinear problem are energy stable provided
\begin{equation}\label{eq:thmEq1}
\mathcal{D}^{+} + {\BCMat}^T\mathcal{D}^{-}{\BCMat} > 0,
\end{equation}
and there exists a positive semi-definite matrix $\Gamma$ such that
\begin{equation}\label{eq:thmEq2}
-\mathcal{D}^{-}+\left(\mathcal{D}^{-}\BCMat\right)\left[\mathcal{D}^{+} + \BCMat^T\mathcal{D}^{-}\BCMat\right]^{-1}\left(\mathcal{D}^{-}\BCMat\right)^T \leq \Gamma < \infty.
\end{equation}
For homogeneous boundary conditions, \eqref{eq:thmEq1} can be relaxed to 
\begin{equation}\label{eq:relaxedBndy}
\mathcal{D}^{+} + {\BCMat}^T\mathcal{D}^{-}{\BCMat} \geq 0,
\end{equation}
and \eqref{eq:thmEq2} becomes redundant.
\end{theorem}
\begin{proof}
See \cite{nordstrom2019}.
\end{proof}

Next, we will discuss the different types of boundaries that can exist.

\subsection{Inflow boundaries}

First, we consider inflow at open boundaries, i.e., $u_n < 0$, $|u_n| > \delta$ {where $\delta$ is a nonzero, possibly small, value}. From the form of the boundary contribution \eqref{eq:rewriteBCs} we require two boundary conditions since
\begin{equation*}
\mathcal{D}^+ = -\frac{\phi}{u_n},
\quad
\mathcal{D}^- =\frac{\phi}{u_n}
\begin{bmatrix}
1 & 0 \\[0.05cm]
0 & 1
\end{bmatrix}.
\end{equation*}
Further, from \eqref{eq:generalFormBCs} we have for the inflow boundary case that
\begin{equation}\label{eq:inflowVars}
\mathbf{w}^-
=
\begin{bmatrix}
u_n^2 + \phi \\[0.05cm]
u_n u_s
\end{bmatrix},
\quad
\BCMat = \begin{bmatrix}
R_1\\[0.05cm]
R_2
\end{bmatrix},
\quad
w^+ = \phi,
\quad
\mathbf{g} = \begin{bmatrix}
g_1\\[0.05cm]
g_2
\end{bmatrix}.
\end{equation}
From the first statement of Theorem~\ref{thm:nonLinStabBCs} we find that
\begin{equation}\label{eq:inflow_first_condition}
\mathcal{D}^{+} + {\BCMat}^T\mathcal{D}^{-}{\BCMat} = \frac{\phi}{u_n}\left(R_1^2 + R_2^2 -1\right)
\end{equation}
becomes a scalar. 
Thus, the coefficients must satisfy the constraint $R_1^2 + R_2^2 < 1$ to guarantee positivity. Additionally, we compute the matrix $\Gamma$ where
\begin{equation*}
\Gamma = -\mathcal{D}^{-}+\left(\mathcal{D}^{-}\BCMat\right)\left[\mathcal{D}^{+} + \BCMat^T\mathcal{D}^{-}\BCMat\right]^{-1}\left(\mathcal{D}^{-}\BCMat\right)^T
=
\frac{\phi}{u_nC}\begin{bmatrix}
R_1^2 - C & R_1R_2\\[0.05cm]
R_1R_2 & R_2^2 - C\\[0.05cm]
\end{bmatrix}
\end{equation*}
where $C = R_1^2 + R_2^2 -1$. To determine if $\Gamma$ is positive semi-definite we compute its eigenvalues to find
\begin{equation*}
\lambda^{\Gamma}_1 = -\frac{\phi}{u_n},\quad\lambda^{\Gamma}_2 = \frac{\phi}{u_n\left(R_1^2+R_2^2-1\right)},
\end{equation*}
which are both positive and bounded because the coefficients for $\BCMat$ are constrained to lie within the unit disc from \eqref{eq:inflow_first_condition} and $|u_n| > \delta$.

\begin{remark}
The analysis above shows that two boundary conditions are sufficient to guarantee a bound on the energy for an inflow boundary. However, it is conceivable that uniqueness might require three. As mentioned above, these issues are presently not clear for general nonlinear hyperbolic problems.
\end{remark}

\begin{remark}
From the linear setting \eqref{eq:diagonalLinearMat}, three boundary conditions are required to bound supercritical flow, i.e. $|\bar{u}_n| > \sqrt{\overline{\phi}}$ {where $ \sqrt{\overline{\phi}} = \sqrt{g\bar{h}}$ is the gravity wave speed}, as
\begin{equation*}
\mathcal{D}^+ = 0 \quad\text{ and }\quad \mathcal{D}^- = \frac{\overline{\phi}}{\bar{u}_n}\begin{bmatrix}
\frac{\bar{u}^2_n - \overline{\phi}}{\overline{\phi}} & 0 & 0\\[0.05cm]
0 & 1  & 0\\[0.05cm]
0 & 0 & 1 \\[0.05cm]
\end{bmatrix}.
\end{equation*}
For subcritical flow, i.e. $|\bar{u}_n| < \sqrt{\overline{\phi}}$, only two are required \cite{huang2014linearized} and the diagonal matrices are
\begin{equation*}
\mathcal{D}^+ = \frac{\overline{\phi}}{\bar{u}_n}\frac{\bar{u}^2_n - \overline{\phi}}{\overline{\phi}} \quad\text{ and }\quad \mathcal{D}^- = \frac{\overline{\phi}}{\bar{u}_n}
\begin{bmatrix}
1  & 0\\[0.05cm]
0 & 1 \\[0.05cm]
\end{bmatrix}.
\end{equation*}
In the nonlinear case, the bound is \textit{not} dependent on the magnitude of the flow and always requires two boundary conditions. From the details given above for the linear case we see that the minimum number of boundary conditions needed for a bound is at least two but possibly three depending on the flow regime. Consequently, the linear and nonlinear analysis are not consistent.
\end{remark}

\subsection{Outflow boundaries}

Next, we examine outflow conditions at open boundaries, i.e., $u_n > \delta > 0$. From the statement \eqref{eq:rewriteBCs} we require one boundary condition since
\begin{equation*}
\mathcal{D}^+ = \frac{\phi}{u_n}
\begin{bmatrix}
1 & 0 \\[0.05cm]
0 & 1
\end{bmatrix},
\quad
\mathcal{D}^- = -\frac{\phi}{u_n},
\end{equation*}
where all entries of the matrix $\mathcal{D}$ are multiplied by the positive quantity $\phi / u_n$.
Thus, the general boundary condition \eqref{eq:generalFormBCs} takes the form
\begin{equation*}
w^- = \phi,
\quad
\BCMat = \begin{bmatrix}
R_1\,\,R_2
\end{bmatrix},
\quad
\mathbf{w}^+
=
\begin{bmatrix}
u_n^2 + \phi \\[0.05cm]
u_n u_s
\end{bmatrix},
\quad
\mathbf{g} = \begin{bmatrix}
g
\end{bmatrix}.
\end{equation*}
To apply the result of Theorem~\ref{thm:nonLinStabBCs} we first compute the matrix
\begin{equation*}
\mathcal{D}^{+} + {\BCMat}^T\mathcal{D}^{-}{\BCMat} = \frac{\phi}{u_n}
\begin{bmatrix}
1-R_1^2 & R_1R_2\\[0.05cm]
R_1R_2 & 1-R_2^2\\[0.05cm]
\end{bmatrix},
\end{equation*}
and calculate its eigenvalues to find
\begin{equation*}
\lambda_1 = \frac{\phi}{u_n},\quad \lambda_2 = \frac{\phi}{u_n}\left(1-R_1^2-R_2^2\right).
\end{equation*}
Here, the constraint $R_1^2 + R_2^2 < 1$ is, again, needed to guarantee positivity. Further, we compute the scalar $\Gamma$ to be
\begin{equation*}
\Gamma = -\mathcal{D}^{-}+\left(\mathcal{D}^{-}\BCMat\right)\left[\mathcal{D}^{+} + \BCMat^T\mathcal{D}^{-}\BCMat\right]^{-1}\left(\mathcal{D}^{-}\BCMat\right)^T =\frac{\phi}{u_n\left(1-R_1^2-R_2^2\right)}
\end{equation*}
which is also positive and bounded whenscaled by $\phi/u_n$ since the coefficients given in $\BCMat$ are contained within the unit disc and $u_n > \delta$.

Alternatively, if we \textit{do not} apply any boundary condition then the boundary term becomes
\begin{equation}\label{eq:useForWallBC}
\mathbf{w}^T\mathcal{D}\,\mathbf{w} = -\frac{\phi^2}{u_n} + \frac{(u_n^2+\phi)^2}{u_n} + u_nu_s^2 = u_n(2\phi + u_n^2 + u_s^2),
\end{equation}
which \textit{still} satisfies the inequality {\eqref{eq:rewriteBCs_inequality}} because of the positive normal velocity. That is, one could view an outflow boundary as a so-called ``free'' boundary \cite{bristeau2001boundary} because no boundary conditions are required.

\begin{remark}
Hence, to obtain an energy bound at outflow boundaries no boundary condition is required. However, uniqueness might still require the imposition of one condition (as was indicated by the eigenvalue analysis above).
\end{remark}

\begin{remark}
In the linear case \eqref{eq:diagonalLinearMat} no boundary condition is needed for supercritical outflow because the diagonal matrices at the boundary are
\begin{equation*}
\mathcal{D}^+ =  \frac{\overline{\phi}}{\bar{u}_n}\begin{bmatrix}
\frac{\bar{u}^2_n - \overline{\phi}}{\overline{\phi}} & 0 & 0\\[0.05cm]
0 & 1 & 0\\[0.05cm]
0 & 0 & 1 \\[0.05cm]
\end{bmatrix}
\quad\text{ and }\quad
\mathcal{D}^- = 0.
\end{equation*}
One boundary condition is needed for subcritical outflow in the linear case \cite{ghader2014} with diagonal matrices
\begin{equation*}
\mathcal{D}^+ =  \frac{\overline{\phi}}{\bar{u}_n}\begin{bmatrix}
1  & 0\\[0.05cm]
0 & 1 \\[0.05cm]
\end{bmatrix}
\quad\text{ and }\quad
\mathcal{D}^- =  \frac{\overline{\phi}}{\bar{u}_n}\frac{\bar{u}^2_n - \overline{\phi}}{\overline{\phi}}.
\end{equation*}

For outflow boundaries, we find, as for inflow boundaries, that the linear and nonlinear analysis are inconsistent. In the nonlinear case, the magnitude of the normal velocity is irrelevant and a bound is obtained with no boundary condition.
\end{remark}

\subsection{Glancing boundaries}

For glancing flow, i.e. $u_n = 0$, there is a singularity in $\mathcal{D}$, however this causes no problem. First, if the flow transitions from outflow to glancing flow, i.e. $u_n\to 0^+$, then 
the formulation \eqref{eq:useForWallBC} shows that a (neutral) energy bound exists for $u_n=0$. Secondly, if the flow transfers from inflow into a glancing flow, i.e. $u_n\to 0^-$, then the formulation \eqref{eq:inflowVars} implies that $R_1 = 1$ and $R_2 = 0$ and the statement \eqref{eq:inflow_first_condition} becomes
\begin{equation*}
\mathcal{D}^{+} + {\BCMat}^T\mathcal{D}^{-}{\BCMat} =  \frac{\phi}{u_n}\left(R_1^2 + R_2^2 -1\right) = 0,
\end{equation*}
which is consistent with \eqref{eq:relaxedBndy}. Hence, Theorem~\ref{thm:nonLinStabBCs} implies that glancing boundaries are (neutrally) stable, and no singularities exist.

\begin{remark}
The nonlinear and linear analysis are consistent for glancing boundaries, i.e. no boundary condition is required provided we know that $u_n=0$. The situation is similar if we instead specify that $u_n=0$ as a boundary condition, the consistency remains. One can intuitively understand this by noting that the additional scalar term in \eqref{eq:nonlinearTerminLemma} vanishes for $u_n=0$. 
\end{remark}

\subsection{Discussion}


We demonstrated that the linear analysis will, in general, yield a different number of boundary conditions compared to the nonlinear problem. In the linear analysis, the number of boundary conditions directly corresponds to the number of negative eigenvalues \cite{bristeau2001boundary,ghader2014,nordstrom2020,shallowwaterbook}. For the linear case this leads to {maximal semi-boundedness and} well-posedness, i.e. a bound on the solution as well as uniqueness and existence. However, in the nonlinear analysis, the {minimal} number of boundary conditions required for an energy estimate could be different than that {determined by} the linear case. In fact, it was possible to obtain an energy estimate by prescribing only two boundary conditions for the inflow case and zero for the outflow case. 

We also found that the magnitude of the normal velocity played \textit{no role} in determining the number of boundary conditions in the nonlinear case. This is quite surprising, again not consistent and quite different from the results of the linear analysis.

Furthermore, we found that the subtleties of the nonlinear energy analysis are not obvious in the traditional entropy analysis. This is because the contraction of the hyperbolic system of conservation laws into entropy space occurs ``too quickly'' and hides information necessary in the boundary condition analysis. We see this directly in the boundary statement from the entropy analysis where the sign of the boundary contribution is dictated entirely by the normal velocity. This provides no clue as to how many conditions one should apply at inflow, in contrast to the number of boundary conditions indicated by the nonlinear energy analysis. 

{From the linear and nonlinear energy analysis in this work we demonstrated that the minimal number of boundary conditions required to guarantee a semi-bounded operator differ. In the linear analysis, maximal semi-boundedness together with certain smoothness assumptions on the initial data is enough to imply well-posedness. For the nonlinear problem maximal semi-boundedness leads to an energy estimate but, unfortunately, it is not known whether it leads to well-posedness since uniqueness and/or existence might be lacking. This is because the underlying well-posedness theory for hyperbolic systems of conservation laws, like the SWEs, is incomplete. In the linear case, the minimal number of boundary conditions together with an energy estimate is \textit{sufficient} to guarantee well-posedness \cite{kreiss1973methods}. For the nonlinear system, maximal semi-boundedness is a \textit{necessary} requirement for well-posedness but it may not be sufficient.}

{Specifically regarding the SWEs, the analysis in this work considered the model without bottom topography terms. For practical simulations in ocean modeling a nonzero bottom topography must be properly discretized and accounted for to develop an energy estimate, see, e.g., \cite{ESDGSEM2D_paper}. This is because the bottom topography contributes to the potential energy in the specific total energy \eqref{eq:totEng}. The analysis in this work is amenable to the case of nonzero bottom topography for the SWEs; however, it requires one to introduce the total geopotential $\Phi=g(h+b)$ as a primitive variable in \eqref{eq:swNoncons} and invoke the assumption that the bottom topography is steady, i.e., $b_t=0$. We plan to extend the nonlinear energy analysis for the SWEs with nonzero bottom topography in future work.}


\section{Concluding remarks}\label{sec:conclusion}

We derived energy stable boundary conditions for the nonlinear shallow water equations. The classical energy method was applied to determine the stability conditions. We demonstrated that there are several possible ways to rewrite terms from the nonlinear analysis into a matrix form. 
Further, we demonstrated that one of these matrix forms was preferable for the boundary condition analysis because it could be rotated into a particularly simple diagonal form. 
It was shown that the traditional energy analysis is consistent with the entropy stability analysis where the  {specific} total energy plays the role of a generalized entropy function. {We also found that the energy analysis provides additional details but is more cumbersome, than the entropy analysis.} 
The details brought forward by the nonlinear energy analysis also made it possible to pinpoint the difference with the linear analysis. 

Specifically, this work investigated and showed:
\begin{enumerate}
   \item The differences and similarities between the nonlinear energy and entropy methods.
   \item Why the linear and nonlinear energy analysis differ and how this affects the boundary conditions.
   \item That a linear analysis cannot be trusted to determine the boundary conditions for a nonlinear problem.
   \item That the flow magnitude does not influence the number of boundary conditions in the nonlinear case.
   \item That the eigenvalues {and their corresponding signs} are not the end all, be all {to indicate the number of boundary conditions necessary for semi-boundedness} in the nonlinear context.
\end{enumerate}

{The use of a nonlinear energy analysis offers two interesting future avenues of study.}
Firstly, for the nonlinear shallow water equations there are interesting theoretical questions concerning {open boundary conditions for flow over} non-constant bottom topography terms as well as possible friction on said bathymetries. This is because the bottom topography contributes to the potential energy of the shallow water equations and the energy analysis requires further investigation. Also, the addition of friction introduces the gradient of the geopotential (or water height) into the model. Therefore, this may expand the possible types of boundary conditions beyond the Dirichlet-type considered herein, to Neumann and Robin type.

Secondly, we want to expand this investigation to other nonlinear hyperbolic systems, like the compressible Euler equations, where the mathematical entropy function is not equal to the  {specific} total energy. {The marriage of traditional energy analysis and nonlinear entropy analysis requires careful attention, particularly when the mathematical entropy function is not the total energy of the system. However, the close relationship between the linear/nonlinear energy analysis and the entropy analysis} could possibly pave the way for truly energy and/or entropy stable approximations. 

\section*{Acknowledgments}

Jan Nordstr\"{o}m was supported by Vetenskapsr{\aa}det, Sweden [award no.~2018-05084 VR] and Andrew Winters by Vetenskapsr{\aa}det, Sweden [award no.~2020-03642 VR].

\bibliographystyle{elsarticle-num}
\bibliography{References}

\appendix
\let\oldsection\section
\renewcommand{\section}[1]{\oldsection{#1}\setcounter{figure}{0}\setcounter{equation}{0}\setcounter{table}{0}}

\section{Proof of Lemma~\ref{lemma:Nmatrix}}\label{sec:newMatrix}
\begin{proof}
To rewrite the scalar term as a quadratic form we assume the existence of a symmetric matrix of the form
\begin{equation*}
\frac{1}{2g}\widetilde{\mathbf{q}}^T\mathcal{N}\widetilde{\mathbf{q}} = \frac{1}{2g}\begin{bmatrix}
\phi\\
u_n\\
u_s
\end{bmatrix}^T
\begin{bmatrix}
n_1 & n_2 & n_3\\[0.05cm]
n_2 & n_4 & n_5\\[0.05cm]
n_3 & n_5 & n_6\\[0.05cm]
\end{bmatrix}
\begin{bmatrix}
\phi\\
u_n\\
u_s
\end{bmatrix}
{=}
\frac{\phi^2}{2g}u_n.
\end{equation*}
The {tangential component of the velocity $u_s$ has} no role in the remaining scalar term (and play no part in the boundary contribution to the energy \cite{nordstrom_roadmap,nordstrom2019}). Thus, $n_3 = n_5 = n_6 = 0$.
Therefore, the matrix ansatz for $\mathcal{N}$ becomes
\begin{equation}\label{eq:matConditionInProof}
\frac{1}{2g}
\begin{bmatrix}
\phi\\
u_n\\
u_s
\end{bmatrix}^T
\begin{bmatrix}
n_1 & n_2 & 0\\[0.05cm]
n_2 & n_4 & 0\\[0.05cm]
0 & 0 & 0\\[0.05cm]
\end{bmatrix}
\begin{bmatrix}
\phi\\
u_n\\
u_s
\end{bmatrix}
=
\frac{1}{2g}\left(n_4 u_n^2 + n_1\phi^2 + 2n_2\phi u_n\right)
{=}
\frac{\phi^2}{2g}u_n.
\end{equation}
This gives three choices on which constant to solve for in the above relationship. First, we solve for $n_1$ to find
\begin{equation*}
n_1 = u_n - \frac{2 n_2 u_n}{\phi} - \frac{n_4 u_n^2}{\phi^2}.
\end{equation*}
This solution yields the matrix
\begin{equation*}
\mathcal{N}_1
=
\begin{bmatrix}
u_n & 0 & 0\\[0.05cm]
0 & 0 & 0\\[0.05cm]
0 & 0 & 0\\[0.05cm]
\end{bmatrix}
+
\begin{bmatrix}
- \frac{2 n_2 u_n}{\phi} - \frac{n_4 u_n^2}{\phi^2} & n_2 & 0\\[0.05cm]
n_2 & n_4 & 0\\[0.05cm]
0 & 0 & 0\\[0.05cm]
\end{bmatrix},
\end{equation*}
where the values of $n_2$ and $n_4$ are arbitrary. However, when multiplying from the left and right with $\widetilde{\mathbf{q}}$ this remaining matrix term vanishes, i.e.,
\begin{equation*}
\frac{1}{2g}\widetilde{\mathbf{q}}^T\mathcal{N}_1\widetilde{\mathbf{q}} = \frac{1}{2g}\widetilde{\mathbf{q}}^T\begin{bmatrix}
u_n & 0 & 0\\[0.05cm]
0 & 0 & 0\\[0.05cm]
0 & 0 & 0\\[0.05cm]
\end{bmatrix}
\widetilde{\mathbf{q}}
+
\frac{1}{2g}
\widetilde{\mathbf{q}}^T
\begin{bmatrix}
- \frac{2 n_2 u_n}{\phi} - \frac{n_4 u_n^2}{\phi^2} & n_2 & 0\\[0.05cm]
n_2 & n_4 & 0\\[0.05cm]
0 & 0 & 0\\[0.05cm]
\end{bmatrix}
\widetilde{\mathbf{q}}
=
\frac{\phi^2}{2g}u_n + 0.
\end{equation*}
Therefore, the values of $n_2$ and $n_4$ are irrelevant to the boundary contribution and we arrive at the first matrix in \eqref{eq:threeMatrices}.
\\
\indent Alternatively, we can elect to solve for $n_2$ in \eqref{eq:matConditionInProof} to obtain
\begin{equation*}
n_2 = \frac{\phi}{2} - \frac{n_1\phi}{2u_n} - \frac{n_4 u_n}{2\phi}.
\end{equation*}
This produces the matrix
\begin{equation*}
\mathcal{N}_2 =  \begin{bmatrix}
0 & \frac{\phi}{2} & 0\\[0.05cm]
\frac{\phi}{2} & 0 & 0\\[0.05cm]
0 & 0 & 0\\[0.05cm]
\end{bmatrix}
+
\begin{bmatrix}
n_1 & - \frac{n_1\phi}{2u_n} - \frac{n_4 u_n}{2\phi} & 0\\[0.05cm]
- \frac{n_1\phi}{2u_n} - \frac{n_4 u_n}{2\phi} & n_4 & 0\\[0.05cm]
0 & 0 & 0\\[0.05cm]
\end{bmatrix},
\end{equation*}
where, now, the values of $n_1$ and $n_4$ remain arbitrary. But when substituting the matrix $\mathcal{N}_2$ into the quadratic form ansatz the remaining matrix vanishes
\begin{equation*}
\frac{1}{2g}\widetilde{\mathbf{q}}^T\mathcal{N}_2\widetilde{\mathbf{q}} = \frac{1}{2g}\widetilde{\mathbf{q}}^T\begin{bmatrix}
0 & \frac{\phi}{2} & 0\\[0.05cm]
\frac{\phi}{2} & 0 & 0\\[0.05cm]
0 & 0 & 0\\[0.05cm]
\end{bmatrix}
\widetilde{\mathbf{q}}
+
\frac{1}{2g}
\widetilde{\mathbf{q}}^T
\begin{bmatrix}
n_1 & - \frac{n_1\phi}{2u_n} - \frac{n_4 u_n}{2\phi} & 0\\[0.05cm]
- \frac{n_1\phi}{2u_n} - \frac{n_4 u_n}{2\phi} & n_4 & 0\\[0.05cm]
0 & 0 & 0\\[0.05cm]
\end{bmatrix}
\widetilde{\mathbf{q}}
=
\frac{\phi^2}{2g}u_n + 0.
\end{equation*}
and the values of $n_1$ and $n_4$ are irrelevant, yielding the second matrix in \eqref{eq:threeMatrices}. 
\\
\indent Lastly, we can solve for the value
\begin{equation*}
n_4 = \frac{\phi^2}{u_n} - \frac{n_1\phi^2}{u_n^2} - \frac{2 n_2 \phi}{u_n},
\end{equation*}
to find the third matrix form
\begin{equation*}
\mathcal{N}_3
=
\begin{bmatrix}
0 & 0 & 0\\[0.05cm]
0 & \frac{\phi^2}{u_n} & 0\\[0.05cm]
0 & 0 & 0\\[0.05cm]
\end{bmatrix}
+
\begin{bmatrix}
n_1 & n_2 & 0\\[0.05cm]
n_2 & - \frac{n_1\phi^2}{u_n^2} - \frac{2 n_2 \phi}{u_n} & 0\\[0.05cm]
0 & 0 & 0\\[0.05cm]
\end{bmatrix}.
\end{equation*}
Substituting the matrix $\mathcal{N}_3$ into the quadratic form, the second matrix term vanishes, i.e.,
\begin{equation*}
\frac{1}{2g}
\widetilde{\mathbf{q}}^T\mathcal{N}_3\widetilde{\mathbf{q}} = \frac{1}{2g}\widetilde{\mathbf{q}}^T\begin{bmatrix}
0 & 0 & 0\\[0.05cm]
0 & \frac{\phi^2}{u_n} & 0\\[0.05cm]
0 & 0 & 0\\[0.05cm]
\end{bmatrix}\widetilde{\mathbf{q}}
+
\frac{1}{2g}
\widetilde{\mathbf{q}}^T
\begin{bmatrix}
n_1 & n_2 & 0\\[0.05cm]
n_2 & - \frac{n_1\phi^2}{u_n^2} - \frac{2 n_2 \phi}{u_n} & 0\\[0.05cm]
0 & 0 & 0\\[0.05cm]
\end{bmatrix}
\widetilde{\mathbf{q}}
=
\frac{\phi^2}{2g}u_n + 0,
\end{equation*}
and we find that the arbitrary values of $n_1$ and $n_2$ are irrelevant. 
\end{proof}

\end{document}